\documentclass[preprint]{elsarticle}

\usepackage{amsthm,amssymb,latexsym,amsxtra}

\usepackage[cmtip,all]{xy}

\newcommand{\lto}{\longrightarrow}

\DeclareMathOperator{\depth}{depth}

\DeclareMathOperator{\Spec}{Spec}
\DeclareMathOperator{\Spm}{Spm}
\DeclareMathOperator{\red}{red}
\DeclareMathOperator{\ann}{ann}
\DeclareMathOperator{\Ass}{Ass}
\DeclareMathOperator{\Supp}{Supp}
\DeclareMathOperator{\modx}{mod}
\DeclareMathOperator{\reg}{reg}

\newtheorem{theorem}{Theorem}
\newtheorem{lemma}[theorem]{Lemma}
\newtheorem{proposition}[theorem]{Proposition}

\theoremstyle{remark}
\newtheorem*{acknowledgements}{Acknowledgements}

\newproof{pf}{Proof}
\newproof{startpf}{Start of proof of Theorem \ref{s1.1}}
\newproof{eproof}{End of proof of Theorem \ref{s1.1}}
\newproof{pot1}{Proof of Theorem \ref{si2}}
\newproof{pot2}{Proof of Theorem \ref{a1.1}}

\begin{document}
\title{A criterion for extending morphisms from open subsets of smooth fibrations of algebraic varieties\tnoteref{t1}}
\tnotetext[t1]{This work was supported by Universit\`{a} di Palermo 
(research project 2012-ATE-0446).\\
Manuscript accepted for publication in Journal of Pure and Applied Algebra,\\ 
https://doi.org/10.1016/j.jpaa.2020.106553 \\
\copyright $\langle 2020\rangle$. This manuscript version is made available under the CC-BY-NC-ND 4.0\\ 
license http://creativecommons.org/licenses/by-nc-nd/4.0/}
\author{Vassil Kanev}
\ead{vassil.kanev@unipa.it}
\date{September 12, 2020}
\address{Dipartimento di Matematica e Informatica,
         Universit\`{a} di Palermo, Via Archirafi 34,
         90123 Palermo, Italy}

\begin{abstract}
Given a smooth morphism $Y\to S$ and a proper morphism $P\to S$ of algebraic varieties we give a sufficient condition for extending 
an $S$-morphism $U\to P$, where $U$ is an open subset of $Y$,
 to an $S$-morphism $Y\to P$, analogous to Zariski's main theorem.
\end{abstract}

\begin{keyword}
extending a morphism \sep Zariski's main theorem \sep smooth fibration \sep fiberwise birational morphism\\
\MSC[2010] 14A10, 14D06
\end{keyword}

\maketitle

\section*{Introduction}
A well-known fact is that every morphism from an open subset of a 
nonsingular algebraic curve to a complete algebraic variety 
may be extended to a morphism of the curve. This is proved by 
considering the closure of the graph of the map in the product of the 
curve and the variety and proving that the projection of the closure 
to the curve is an isomorphism by means of Zariski's main theorem. 
It is desirable to generalize this extension property to families of curves,
with a view toward applications to moduli problems. 
Trying to apply the same argument to families of curves one 
encounters the following problem. Given two families of curves 
$g:X\to S$, $h:Y\to S$, where $h$ is smooth, and a commutative diagram of 
morphisms as in \eqref{ei2} below,
such that $f$ is finite and surjective, and 
$f_{s}:g^{-1}(s)\to h^{-1}(s)$ is birational for every $s\in S$, is it 
true that $f$ is an isomorphism? If $g$ and $h$ were proper and if the scheme-theoretic fibers of $g$ were 
reduced this would follow from Proposition~4.6.7~(i) of 
\cite{EGAIII.1}. The latter condition, however, might be difficult to 
verify. In fact a weaker condition on $g$ suffices to conclude that 
$f$ is an isomorphism. Before we formulate one of our results let us recall 
that a morphism of algebraic varieties $h: Y\to S$ over an algebraically closed field $k$
is smooth at $y\in Y$ if $h$ is flat at $y$ and
$\dim_{k} \Omega_{Y/S}(y)\leq \dim_{y}Y_{h(y)}$.  The following theorem holds.
\renewcommand{\theenumi}{\alph{enumi}}
\begin{theorem}\label{si2}
Let $X,Y$ and $S$ be algebraic varieties over an algebraically closed field $k$.
Let $d$ be an integer $\geq 1$. Let
\begin{equation}\label{ei2}
\xymatrix{
X\ar[rr]^{f}\ar[rd]_{g}&&
Y\ar[dl]^{h}\\
&S
}
\end{equation}
be a commutative diagram of morphisms such that: 
\begin{enumerate}
\item
$h: Y\to S$ is smooth of relative dimension $d$;
\item $f$ is finite and surjective;
\item
for every $s\in h(Y)$ the fiber $g^{-1}(s)$ is irreducible and there 
is  point $x\in g^{-1}(s)$ such that $g$ is smooth at $x$;
\item
for every $s\in g(X) = h(Y)$ the map $f_{s}:g^{-1}(s)\to h^{-1}(s)$ 
is birational.
\end{enumerate}
Then $f:X\to Y$ is an isomorphism. 
\end{theorem}
 Theorem~\ref{si2} is used
in the proof of the following criterion for extension of morphisms. 
\begin{theorem}\label{a1.1}
Let $Y, S$ and $P$ be algebraic varieties over an algebraically closed field 
$k$. Let $d$ be an integer $\geq 1$. Let $h:Y\to S$  be  a  smooth  morphism 
whose nonempty fibers are irreducible of dimension $d$. Let $P\to  S$  be  a 
proper morphism. Let $U$ be an open subset of $Y$ such that 
$U\cap h^{-1}(s)\neq \emptyset$ for every $s\in h(Y)$.  Let  $\varphi  :U\to 
P$ be an $S$-morphism. Let $\Gamma \subset  U\times_{S}P$  be  its graph
and let $X = \overline{\Gamma} \subset  Y\times_{S}P$ be its closure. 
Suppose that the projection $f:X\to Y$ has the property that 
$f^{-1}(y)$ is a finite set for every $y\in Y\setminus U$. Then there  is  a 
unique extension of $\varphi$ to $Y$: an $S$-morphism 
$\tilde{\varphi}:Y\to P$ such that $\tilde{\varphi}|_{U}=\varphi$.
\end{theorem}
We notice that if $S$ is irreducible and normal, then $Y$ is irreducible and normal as well, and the statement of Theorem~\ref{a1.1} 
follows from Zariski's main theorem.
\par
Under the assumptions of Theorem~\ref{si2}, when the relative dimension of $g: X\to S$ is one,  
Theorem~9  of  \cite{Ko} 
applied to $\Spec\mathcal{O}_{X,x} \to \Spec\mathcal{O}_{S,s}$, where 
$x\in X$ is arbitrary, $s=g(x)$, yields that $g:X\to S$ is smooth, which
implies that 
$f: X\to Y$ is an isomorphism. The proof of Theorem~9  of  \cite{Ko} needs however
an additional assumption, as the author was kindly informed by Prof. J. Koll\'{a}r (letter of January 21, 2019), 
namely Assumption~9.2 of Theorem~9 (ibid.) should hold for the completions 
$(\Spec\mathcal{O}_{X,x})\sphat \to (\Spec\mathcal{O}_{S,s})\sphat$. This condition might
not hold under the assumptions of Theorem~\ref{si2}.
\par
The paper is organized as follows. In Section~\ref{s1} we prove Theorem~\ref{s1.1}, which is a 
result of  commutative 
algebra  of  independent  interest,  and  from  which   Theorem~\ref{si2}   and 
Theorem~\ref{a1.1} are deduced. Given a local flat homomorphism  of  Noetherian 
local rings $(A,\mathfrak{m})\to (R,\mathfrak{n})$ with regular fiber 
$R/\mathfrak{m}R$ and a finite homomorphism $R\to B$, 
Theorem~\ref{s1.1} gives sufficient conditions which ensure that $R\to B$ is  an 
isomorphism. The proof of this theorem is  based on the key idea of  J. Koll\'{a}r
by which Theorem~9 of \cite{Ko} is proved, namely splitting the $R$-module $B$ into
a direct sum of $R$ and another summand, then proving that this summand is zero.
In Section~\ref{s2} we give the proofs of Theorem~\ref{si2} and Theorem~\ref{a1.1}. 

\medskip \noindent
\emph{Notation.} If $\varphi: A\to B$ is a homomorphism of rings and $I\subset A$, $J\subset B$ 
are ideals we denote, following  \cite{Ma2},  $\varphi  (I)B$  by  $IB$  and 
$\varphi^{-1}J$ by $A\cap J$. If $\mathfrak{p}\subset A$ is a  prime  ideal, 
then $k(\mathfrak{p})$ is the quotient field of $A/\mathfrak{p}$.
A variety over an algebraically closed field $k$ is a reduced separated, possibly reducible, scheme 
of finite type over $k$. A point of a variety is a closed point of the scheme. The maximal spectrum 
of a ring $A$ is denoted by $\Spm A$.

\section{The main theorem}\label{s1}
\begin{theorem}\label{s1.1}
Let $(A,\mathfrak{m})$ and $(R,\mathfrak{n})$ be Noetherian  local  rings, and let 
 $A\to R$ be a local flat  homomorphism.  Suppose that $A$  is  reduced. 
Suppose that:
\renewcommand{\theenumi}{\alph{enumi}}
\begin{enumerate}
\item
$R\otimes_{A}k(\mathfrak{m})$ is a regular  ring of dimension $n\geq 1$ and  $k(\mathfrak{n})$  is 
separable over $k(\mathfrak{m})$;
\item
$R\otimes_{A}k(\mathfrak{p})$ is a regular  ring  for  every  minimal  prime 
ideal $\mathfrak{p}$ of $A$.
\end{enumerate}
Set $S=\Spec A$, $s_{0}=\mathfrak{m}$, $Y=\Spec R$,  $y_{0}=\mathfrak{n}$, 
and let $h:Y\to S$ be the associated morphism  of  affine  schemes.  Suppose 
there is a commutative diagram of morphisms of schemes 
\begin{equation}\label{es1.1}
\xymatrix{
X\ar[rr]^{f}\ar[rd]_{g}&&
Y\ar[dl]^{h}\\
&S
}
\end{equation}
such that:
\renewcommand{\theenumi}{\roman{enumi}}
\begin{enumerate}
\item
$f$ is finite and surjective;
\item
$X_{s_0}$ is irreducible and generically reduced;
\item
$g$ is flat at the generic point of $X_{s_{0}}$;
\item
if  $\eta  \in  X_{s_{0}}$  and  $\zeta  =  f(\eta)   \in   Y_{s_{0}}$   are 
the    generic    points    of    $X_{s_{0}}$    and    $Y_{s_{0}}$     then 
\linebreak
$f^{\sharp}(\zeta):k(\zeta)\to k(\eta)$ is an isomorphism;
\item
every irreducible component of $X$ contains $g^{-1}(s_{0})$.
\end{enumerate}
Then $Y$ is reduced and $f\circ i:X_{\red}\to Y$ is an isomorphism.  Moreover 
the open set of  reduced  points  of  $X$  contains  the  generic  point  of 
$g^{-1}(s_{0})$.  Assumption~(a)  implies  Assumption~(b)  if   $R$   is   a 
localization of a finitely generated $A$-algebra, or if $A$  is  a  $G$-ring 
(cf. \cite[\S~34]{Ma1}, or \cite[\S~32]{Ma2}).
\end{theorem}
\begin{startpf}
Let $\mathfrak{p}_{1},\ldots,\mathfrak{p}_{m}$ be the minimal  prime  ideals 
of $A$. One has  $A_{\mathfrak{p}_{i}}=k(\mathfrak{p}_{i})$  for  every  $i$ 
since $A$ is reduced. Tensoring $0\to A\to  \oplus_{i}A_{\mathfrak{p}_{i}}$ 
by $R$ one obtains by Assumption~(b) that  $R$  is  reduced.  Assumption~(i) 
implies that $X\cong \Spec B$, where  $B=\mathcal{O}_{X}(X)$  is  a  finite 
$R$-module and moreover $f^{\sharp}(Y):R\to B$ is injective,  since  $f$  is 
surjective and $R$ is reduced. The points of $X$ where $f$ is not flat  form 
a closed subset $Z\subset X$. The image $f(Z)$ is closed in $Y$ since $f$ is 
finite. Let $Y'=Y\setminus f(Z)$, $X'=f^{-1}(Y')$. Then  $f|_{X'}:X'\to  Y'$ 
is  finite  and  flat.  
\par
We  claim  that  $\zeta   \in   Y'$.   One   applies 
\cite[20.G]{Ma1}   to   $\mathcal{O}_{S,s_{0}}\to   \mathcal{O}_{Y,\zeta}\to 
\mathcal{O}_{X,\eta}$   and    $M=\mathcal{O}_{X,\eta}$.    By    hypothesis 
$\mathcal{O}_{Y,\zeta}$     and     $\mathcal{O}_{X,\eta}$     are      flat 
$\mathcal{O}_{S,s_{0}}$-modules.  Furthermore  $Y_{s_{0}}$  is  integral   by 
Assumption~(a),               so 
$\mathcal{O}_{Y,\zeta}\otimes _{\mathcal{O}_{S,s_{0}}}k(s_{0})\cong 
\mathcal{O}_{Y_{s_{0}},\zeta}\cong 
k(\zeta)$             is             a             field.              Hence 
$\mathcal{O}_{X,\eta}\otimes_{\mathcal{O}_{S,s_{0}}}k(s_{0})$ is a flat 
$\mathcal{O}_{Y,\zeta}\otimes      _{\mathcal{O}_{S,s_{0}}}k(s_{0})$-module. 
Therefore $\mathcal{O}_{X,\eta}$ is a  flat  $\mathcal{O}_{Y,\zeta}$-module. 
The  hypothesis   that   $g^{-1}(s_{0})$   is   irreducible   implies   that 
$f^{-1}(\zeta)=\eta$, therefore $\zeta \in Y', \eta \in X'$.
\par
Assumption~(i) and Assumption~(v) imply that every irreducible 
component of $Y$, being an image of some irreducible component of $X$, 
contains $h^{-1}(s_{0})$ and in particular $\zeta$. Therefore the open 
set $Y'$ is connected and dense in $Y$. We claim that $f|_{X'}:X'\to 
Y'$ is an isomorphism. First, $f_{*}\mathcal{O}_{X'}$ is a locally 
free sheaf of a certain rank $d\geq 1$, since $f|_{X'}$ is finite, 
surjective and flat, and $Y'$ is connected. One has 
$d=\dim_{k(y)}\Gamma(X_{y},\mathcal{O}_{X_{y}})$ for every $y\in Y'$. 
Let $y=\zeta$. Then $f^{-1}(y)=\eta$ and furthermore $f:X\to Y$ is 
unramified at $\eta$. Indeed, it suffices to verify that 
$f|_{X_{s_{0}}}:X_{s_{0}}\to Y_{s_{0}}$ is unramified at $\eta$. This 
holds since $\mathcal{O}_{Y_{s_{0}},\zeta} = k(\zeta)$, 
$\mathcal{O}_{X_{s_{0}},\eta} = k(\eta)$, for by hypothesis 
$X_{s_{0}}$ is irreducible and generically reduced, and furthermore 
$k(\zeta)\to 
k(\eta)$ is an isomorphism by Assumption~(iv). We obtain that for $y=\zeta$, 
$X_{y}=\Spec k(\eta)$ and 
$d=\dim_{k(y)}\Gamma(X_{y},\mathcal{O}_{X_{y}})=1$. This shows that the 
morphism of sheaves of rings 
$f^{\sharp}:\mathcal{O}_{Y'}\to f_{*}\mathcal{O}_{X'}$ makes 
$f_{*}\mathcal{O}_{X'}$ a locally free 
$\mathcal{O}_{Y'}$-module of rank 1. This implies that 
$f^{\sharp}:\mathcal{O}_{Y'}\to f_{*}\mathcal{O}_{X'}$ is an 
isomorphism, hence $f|_{X'}:X'\to Y'$ is an isomorphism. Since $Y$ is 
reduced this implies that the open set of reduced points of $X$ 
contains $X'$, in particular the generic point $\eta$ of $g^{-1}(s_{0})$ is 
a reduced point of $X$.
\par
In order to prove the isomorphism 
$f\circ i:X_{\red}\overset{\sim}{\lto} Y$ we replace $X$ by $X_{\red}$ and 
observe that all the assumptions of the theorem hold for 
$f\circ i = f_{\red}:X_{\red}\to Y$ and $g_{\red}:X_{\red}\to S$. We may thus 
assume 
that $X=\Spec B$ is reduced, so by Assumption~(v)  
$g^{-1}(s_{0})\subset V(\mathfrak{p})$, where $\mathfrak{p}$ is any associated prime ideal of $B$.
\par
Let $I$ be the radical ideal $I=I(Y\setminus Y')\subset R$. We will 
prove below  that if $I\ne R$, then $\depth (I,R)\geq 2$. Assuming this statement one proves that $f:X\to Y$ is an 
isomorphism as follows. Consider the exact sequence of finite 
$R$-modules
\begin{equation}\label{es1.3}
0 \to R \overset{f^{\sharp}(Y)}{\lto} B \to Q \to 0.
\end{equation}
Let $R'$ be the image of $R$. By way of contradiction let us assume 
$Q\ne 0$. Let $J=\surd \ann(Q)$. 
One has $V(J)=\Supp\, Q$, so $Y\setminus V(J)\subset Y'$. The 
isomorphism $f^{-1}(Y')\overset{\sim}{\lto} Y'$, proved above, shows 
that $Y' \subset Y\setminus V(J)$. Therefore $I=J\ne R$. The stated 
inequality $\depth(I,R)\geq 2$ implies by \cite[Theorem~16.6]{Ma2} that 
$Ext^{1}_{R}(Q,R)=0$. Hence \eqref{es1.3} splits, $B\cong R'\oplus 
Q'$. Let $\mathfrak{p}=\ann_{R}(x)$ be an associated prime ideal of 
$Q$. Since $\Ass(Q)\subset \Supp(Q) = V(J)$ (cf. 
\cite[Theorem~6.5]{Ma2}), one has $\mathfrak{p}\supset I$. Let 
$\mathfrak{p}'\subset R'$ be the image of $\mathfrak{p}$ and let 
$x'\in Q'$ be the preimage of $x$. Since $\mathfrak{p}'\cdot x'=0$ the 
subset $\mathfrak{p}'\subset B$ consists of zero divisors, so 
$\mathfrak{p}'\subset P_1\cup 
\cdots \cup P_{m}$, where $P_{i}, i=1,\ldots,m$, are the associated prime 
ideals of $B$. Let $\mathfrak{p}_{i}= R\cap P_{i}$. Then 
$\mathfrak{p}\subset \mathfrak{p}_{1}\cup \cdots \cup 
\mathfrak{p}_{m}$, so $\mathfrak{p}\subset \mathfrak{p}_{j}$ for some 
$j$. Let $P=\surd(\mathfrak{m}B)$. Assumption~(v) means that 
$P_{i}\subset P$ for every $i$. Hence $R\cap P\supset 
\mathfrak{p}_{j}\supset \mathfrak{p} \supset I$. Since $P\in \Spec B$ 
and $R\cap P \in \Spec R$ are the same as $\eta \in X$ and $\zeta \in 
Y$ respectively, one obtains that $\zeta$ belongs to $V(I)=Y\setminus 
Y'$. This contradiction shows that $Q=0$ and therefore $f^{\sharp}(Y): 
R\to B$ is an isomorphism. The isomorphism $f\circ i:X_{\red}\to Y$ is 
proved.
\par
We prove now the last statement of the theorem. Let $C$ be a finitely 
generated  $A$-algebra. Let $T=\Spec C$, $u:T\to S=\Spec A$ be the morphism 
corresponding to $A\to C$. Suppose there is an $A$-isomorphism  of  $R$  with 
$\mathcal{O}_{T,z}$ for some $z\in T$. Let $j:Y\to  T$  be  the  composition 
$Y\overset{\sim}{\lto} \Spec \mathcal{O}_{T,z}\lto T$. One has $u(z)=s_{0}$ 
and $\mathcal{O}_{T,z}\otimes_{\mathcal{O}_{S,s_{0}}}k(s_{0})\cong 
\mathcal{O}_{T_{s_{0}},z}$. Assumption~(a) implies, according to 
\cite[Corollaire~II.5.10]{SGA1}, 
that the  fiber  $T_{s_{0}}$  is  smooth  at  $z$, so $u$ is smooth at $z$
by \cite[Theorem~VII.1.8]{A-K}. 
Smoothness at a point is an open condition so  we  may,  replacing  $T$  by  an  affine 
neighborhood of $z$, assume that  $u:T\to  S$  is  smooth.  It  is  moreover 
surjective since the flat morphism $u$ is an open map. Every  fiber  $T_{s}, 
s\in S$ is geometrically regular. Let $s\in S, y\in h^{-1}(s), t=j(y)$.  One 
has $\mathcal{O}_{S,s}$-isomorphisms
\[
\mathcal{O}_{Y_{s},y} \cong \mathcal{O}_{Y,y}\otimes_{\mathcal{O}_{S,s}}k(s) 
\cong \mathcal{O}_{T_{s},t}.
\]
Therefore $\mathcal{O}_{Y_{s},y}$ is geometrically regular. This means  that 
$R\otimes_{A}k(\mathfrak{p})$  is  geometrically  regular  ring  for   every 
$\mathfrak{p}\in \Spec A$, in particular Assumption~(b) holds.
\par
Suppose now that $A$  is  a  $G$-ring.  Then  $A$  is  quasi-excellent  (cf. 
\cite[\S~34]{Ma1}). Let $k=k(\mathfrak{m})$, 
$\mathfrak{n}_{0}=\mathfrak{n}/\mathfrak{m}R\subset R\otimes_{A}k$. 
Assumption~(a) implies that $R\otimes_{A} k$ is formally smooth with  respect 
to the $\mathfrak{n}_{0}$-adic topology  (cf.  
\cite[\S~28.M]{Ma1}~Proposition). 
Hence by Th\'{e}or\`{e}m~19.7.1  of  \cite[Ch.0]{EGAIV.1}  the  homomorphism 
$A\to R$ is formally smooth with  respect  to  the  $\mathfrak{m}$-adic  and 
$\mathfrak{n}$-adic topologies of  $A$  and  $R$.  A  theorem  of  Andr\'{e} 
\cite{Andre}  yields  that  $A\to  R$  is   a   regular   homomorphism,   so 
$R\otimes_{A}k(\mathfrak{p})$   is   geometrically   regular    for    every 
$\mathfrak{p}\in \Spec A$. This implies Assumption~(b).
\end{startpf}
Our next goal is to prove that $\depth (I,R)\geq  2$  provided  $I\ne  R$,  a 
statement  used  in  the  proof  of  Theorem~\ref{s1.1}.  It  is  proved  in 
\cite[Ch.0]{EGAIII.1} Proposition~10.3.1 that, given a Noetherian local  ring 
$(A,\mathfrak{m})$ and a  homomorphism  of  fields  $k(\mathfrak{m})\to  K$, 
there exists a Noetherian local ring $(B,J)$ and a flat  local  homomorphism 
$(A,\mathfrak{m})\to (B,J)$  such  that  $\mathfrak{m}B=J$  and  $k(J)$  is 
isomorphic  to  $K$  over  $k(\mathfrak{m})$.  We   include  the 
proof of the following known fact since we could not find a reference.
\begin{lemma}\label{s1.6}
Let $(A,\mathfrak{m})$ and $(R,\mathfrak{n})$ be Noetherian  local  rings, and let 
 $A\to R$ be a local flat  homomorphism.  Let   $k=k(\mathfrak{m})$,   $K=k(\mathfrak{n})$. 
Suppose $K$ is separable over $k$.  Suppose  $R\otimes_{A}k$  is  a  regular 
ring of dimension $n\geq 1$. Let $(B,J)$ and $A \to B$ be as above: the homomorphism is local and  flat, 
$\mathfrak{m}B=J$ and $k(J)$  is  $k$-isomorphic  to  $K$.  Then  there  is  an 
$A$-isomorphism    $\hat{R}\cong    \hat{B}[[T_{1},\ldots,T_{n}]]$,    where 
$\hat{R}$ is the $\mathfrak{n}$-adic completion of $R$ and $\hat{B}$ is  the 
$J$-adic completion of $B$.
\end{lemma}
\begin{pf}
Let us first consider the case where $k(\mathfrak{m})\to k(\mathfrak{n})$  is 
an isomorphism and $(B,J)=(A,\mathfrak{m})$. We have a  commutative  diagram 
of faithfully flat homomorphisms (cf. \cite[Theorem~22.4]{Ma2})
\begin{equation}\label{es1.6}
\xymatrix{
A\ar[r]\ar[d]&R\ar[d]\\
\hat{A}\ar[r]&\hat{R}
}
\end{equation}
Let   $t_{1},\ldots,t_{n}\in   \mathfrak{n}$   be   elements    such    that 
$x_{i}=t_{i}(\modx\mathfrak{m}R), i=1,\ldots,n$ generate the maximal ideal  
of 
the   regular   local   ring   $R/\mathfrak{m}R\cong   R\otimes_{A}k$.   Let 
$\varphi :\hat{A}[[T_{1},\ldots,T_{n}]]\lto \hat{R}$ be the homomorphism  of 
$\hat{A}$-algebras      such      that      $\varphi(T_{i})=t_{i}$      (cf. 
\cite[Theorem~7.16]{Ei}). Let  $\hat{\mathfrak{m}}=\mathfrak{m}\hat{A}$  and 
$\hat{\mathfrak{n}}=\mathfrak{n}\hat{R}$ be the maximal ideals of  $\hat{A}$ 
and $\hat{R}$. The  ring  $A'=\hat{A}[[T_{1},\ldots,T_{n}]]$  is  local  and 
complete with maximal ideal $M=(\hat{\mathfrak{m}},T_{1},\ldots,T_{n})$ (cf. 
\cite[Ch. III, \S~2.6, Proposition~6]{Bou}). One  has  $\varphi(M)=\hat{\mathfrak{n}}$,  so 
$\hat{R}/M\hat{R}\cong \hat{R}/\hat{\mathfrak{n}}\cong R/\mathfrak{n}  \cong 
k(\mathfrak{n}) \cong k$. Hence $\varphi:A'\to \hat{R}$ is  surjective  (cf. 
\cite[Theorem~8.4]{Ma2}). In order to  prove  that  $\varphi$  is  injective, 
applying    \cite[Theorem~22.5]{Ma2},    we    need    to    verify     that 
$\overline{\varphi}:A'\otimes_{\hat{A}}k\to  \hat{R}\otimes_{\hat{A}}k$   is 
injective.         One          has          $\hat{R}\otimes_{\hat{A}}k\cong 
\hat{R}/\hat{\mathfrak{m}}\hat{R}\cong    \widehat{R/\mathfrak{m}R}$.     By 
assumption   the   composition   $k=k(\mathfrak{m})\to    R/\mathfrak{m}R\to 
R/\mathfrak{n}=k(\mathfrak{n})$ is an isomorphism. Hence 
$\widehat{R/\mathfrak{m}R}\cong k[[x_{1},\ldots,x_{n}]]$ (cf. \cite[p.~124, 
Remark~2]{A-M}). Therefore the composition 
\begin{equation*}
k[[T_{1},\ldots,T_{n}]]     \overset{\sim}{\lto}     A'/\hat{\mathfrak{m}}A' 
\overset{\overline{\varphi}}{\lto}         \hat{R}/\hat{\mathfrak{m}}\hat{R} 
\overset{\sim}{\lto} k[[x_{1},\ldots,x_{n}]]
\end{equation*}
which transforms $T_{i}$ in $x_{i}$ is an  isomorphism.  This  implies  that 
$\overline{\varphi}: A'\otimes_{\hat{A}}k\to  \hat{R}\otimes_{\hat{A}}k$  is 
an isomorphism.  We conclude that $\varphi :\hat{A}[[T_{1},\ldots,T_{n}]]\to 
\hat{R}$ is an isomorphism provided $k(\mathfrak{m})\to k(\mathfrak{n})$  is 
an isomorphism.
\par
Let  us  consider  now   the   general   case   when   $k=k(\mathfrak{m})\to 
k(\mathfrak{n})=K$    is    an    arbitrary    separable    extension.    By 
\cite[Theorem~26.9]{Ma2}   $K$   is   $0$-smooth   over   $k$.   Hence    by 
\cite[Theorem~28.10]{Ma2} $B$ is $J$-smooth over $A$. This implies that  the 
homomorphism $B\to K=\hat{R}/\hat{\mathfrak{n}}$ has a lifting $\varphi:B\to 
\hat{R}$   which   is   a   local   homomorphism   of   $A$-algebras    (see 
\cite[p.~214]{Ma2}).  Applying  \cite[20.G]{Ma1}  to  $A\to  B\to  \hat{R}$, 
taking into account that  $B\otimes_{A}k\to  \hat{R}\otimes_{A}k$  is  flat 
since $B\otimes_{A}k\cong K$ is a field, we conclude that $B\to \hat{R}$  is 
flat.     Furthermore      $k(J)=B/J\to      \hat{R}/\hat{\mathfrak{n}}\cong 
R/\mathfrak{n}= k(\mathfrak{n})$ is an isomorphism and by 
\cite[Theorem~8.11]{Ma2}
\begin{equation*}
\hat{R}\otimes_{B}k(J) \cong \hat{R}\otimes_{B}(B\otimes_{A}k) \cong 
\hat{R}\otimes_{A}k \cong \hat{R}/\mathfrak{m}\hat{R} \cong 
\widehat{R/\mathfrak{m}R}
\end{equation*}
is a regular local ring of dimension $n$. By the first part of the proof one 
concludes that $\hat{R}\cong \hat{B}[[T_{1},\ldots,T_{n}]]$
\end{pf}
\begin{proposition}\label{s1.8}
Let $(A,\mathfrak{m})$ and $(R,\mathfrak{n})$ be Noetherian  local  rings, and let 
 $A\to R$ be a local flat  homomorphism.
Suppose $k=k(\mathfrak{m})\to k(\mathfrak{n})= K$ is a separable extension. 
Suppose $R\otimes_{A}k$ is a regular ring. Let $I\subset R$ be a proper 
ideal such that:
\renewcommand{\theenumi}{\alph{enumi}}
\begin{enumerate}
\item
none of the prime ideals of the set $A\cap V(I)\subset \Spec A$ is 
contained in an associated prime ideal of $A$;
\item
$I \not \subset \mathfrak{m}R$.
\end{enumerate}
Then $\depth(I,R)\geq 2$.
\end{proposition}
\begin{pf}
We may replace $I$ by its radical and thus assume that  $I=\surd I$.  Indeed, 
$\depth(I,R)=\depth(\surd I,R)$         (see         \cite[Corollary~17.8]{Ei}), 
$V(I)=V(\surd I)$ and $I\not \subset \mathfrak{m}R$ if and only if 
$\surd I\not \subset \mathfrak{m}R$ since the condition that  $R\otimes_{A}k$ 
is regular implies that $\mathfrak{m}R$ is a prime ideal. Furthermore Condition~(b)
implies that $\mathfrak{m}R \subsetneqq \mathfrak{n}$, so $\dim R\otimes_{A}k \geq 1$.
 Let  $I=P_{1}\cap 
\cdots \cap P_{r}$, where $P_{i}, i=1,\ldots,r$ are the minimal prime ideals 
which contain $I$. We claim that there exists $a_{1}\in A\cap I$, such  that 
$a_{1}$ is not a zero divisor of $A$. If this were not the case, then $A\cap 
I\subset   \mathfrak{p}_{1}\cup   \cdots   \cup   \mathfrak{p}_{s}$,   where 
$\mathfrak{p}_{i}, i=1,\ldots,s$ are the  associated  primes  of  $A$.  Then 
$A\cap I\subset \mathfrak{p}_{j}$  for  some  $j$  and  consequently  $A\cap 
P_{i}\subset  \mathfrak{p}_{j}$  for  some  $i$ (cf. \cite[Proposition~1.11]{A-M}).  
This  contradicts Condition~(a). 
\par
Let $\hat{I}$,  $\hat{R}$  be  the  $\mathfrak{n}$-adic 
completions of $I$ and $R$. The equality $\depth(I,R)=\depth(\hat{I},\hat{R})$ 
holds. Indeed, let $I=(x_{1},\ldots,x_{m})$.  Consider  the  Koszul  complex 
$K^{\bullet}=K^{\bullet}(x_{1},\ldots,x_{m})$,                               
where 
$K^{i}(x_{1},\ldots,x_{m})=\Lambda^{i}N,    N=\oplus_{i=1}^{m}Re_{i}$    and 
$d^{i}:K^{i}\to     K^{i+1}$     is     $d^{i}(v)=x\wedge     v$,      where 
$x=\sum_{i=1}^{m}x_{i}e_{i}$. Then  $\depth(I,R)=r$  iff  
$H^{i}(K^{\bullet})=0$  for  $i<r$  and   $H^{r}(K^{\bullet})\ne   0$   (cf. 
\cite[Theorem~17.4]{Ei}). Since $\hat{I}=I\hat{R}$ the  images  $x'_{i}$  of 
$x_{i}$ in $\hat{R}$, $i=1,\ldots,m$, generate $\hat{I}$.  The  corresponding 
Koszul      complex      is      $K^{'\bullet}=K'(x'_{1},\ldots,x'_{m})\cong 
K^{\bullet}(x_{1},\ldots,x_{m})\otimes_{R}\hat{R}$. Since $R\to \hat{R}$  is 
faithfully     flat     one     has      that      $H^{i}(K^{'\bullet})\cong 
H^{i}(K^{\bullet})\otimes_{R}\hat{R}$ and $H^{i}(K^{'\bullet}) \neq 0$ if and only if 
$H^{i}(K^{\bullet}) \neq 0$.    
By     the     above     criterion 
$\depth(I,R)=\depth(\hat{I},\hat{R})$. 
\par
Let       $(A,\mathfrak{m})\to        (B,J)$        and        $\hat{R}\cong 
\hat{B}[[T_{1},\ldots,T_{n}]]$ be  as  in  Lemma~\ref{s1.6}.  The  hypothesis 
$I\not  \subset   \mathfrak{m}R$   implies   $\hat{I}=I\hat{R}\not   \subset 
(\mathfrak{m}R)\hat{R}=\mathfrak{m}\hat{R}$ since $R\to  \hat{R}$,  being  a 
faithfully flat homomorphism, has the property that $\mathfrak{a}\hat{R}\cap 
R = \mathfrak{a}$ for every ideal $\mathfrak{a}$ in $R$. Furthermore the ring 
extensions      $A\to       B\to       \hat{B}\to       \hat{R}$       yield 
$\mathfrak{m}\hat{R}=(\mathfrak{m}B)\hat{R} =  J\hat{R}  =  \hat{J}\hat{R}$. 
Therefore $\hat{I}\not \subset \hat{J}\hat{R}$. 
\par
Let $a_{1}\in A\cap I$ be  a non zero divisor of $A$ as above. Let 
$a'_{1}\in \hat{I}$ be  its  image  in 
$\hat{R}$.  Let  $f\in  \hat{I}\setminus  \hat{J}\hat{R}$.  We  claim   that 
$a'_{1},f$   is   an   $\hat{R}$-regular   sequence.   This   implies   that 
$\depth(\hat{I},\hat{R})\geq 2$. Abusing notation we identify $\hat{R}$  with 
$\hat{B}[[T_{1},\ldots,T_{n}]]$.          Let          $\overline{f}=f(\mod\hat{J}\hat{R})\in 
K[[T_{1},\ldots,T_{n}]]$. There exist  positive  integers 
$u_{1},u_{2},\ldots,u_{n}$   such   that    the    automorphism    $s$    of 
$K[[T_{1},\ldots,T_{n}]]$  defined  by  $s(T_{i})=T_{i}+T_{n}^{u_{i}}$   for 
$1\leq  i\leq  n-1$  and  $s(T_{n})=T_{n}$  transforms   $\overline{f}$   in 
$\overline{g}(T_{1},\ldots,T_{n})$  with  $\overline{g}(0,\ldots,0,T_{n})\ne 
0$ (cf. \cite[Ch.VII \S~3 no.7]{Bou} Lemma~3). The same substitution yields  a 
$\hat{B}$-automorphism  $\varphi$  of  $\hat{B}[[T_{1},\ldots,T_{n}]]$.  Let 
$g=\varphi(f)$. Let $C=\hat{B}[[T_{1},\ldots,T_{n-1}]]$. This is a  complete 
local ring with maximal ideal  $\mathfrak{M}=(\hat{J},T_{1},\ldots,T_{n-1})$ 
and $\hat{R}\cong  C[[T_{n}]]$.  We  claim  that  $a'_{1},g$  is  a  regular 
$\hat{R}$-sequence.   Indeed,   the    injectivity    of 
$A\overset{\cdot a_{1}}{\lto}A$ implies the injectivity of 
$\hat{R}\overset{\cdot a_{1}}{\lto}\hat{R}$ since the composition $A\to R\to 
\hat{R}$ is flat. One has          
$g(\modx\mathfrak{M})=\overline{g}(0,\ldots,0,T_{n})\ne 0$     and 
$\hat{R}/a'_{1}\hat{R}\cong (\hat{B}/a'_{1}\hat{B})[[T_{1},\ldots,T_{n}]]
\cong \overline{C}[[T_{n}]]$, where $\overline{C}=
(\hat{B}/a'_{1}\hat{B})[[T_{1},\ldots,T_{n-1}]]$ is a complete local ring 
with maximal   ideal   $(\hat{J}/a'_{1}\hat{J},T_{1},\ldots,T_{n-1})$.   
Applying 
\cite[Ch.VII \S~3 no.8]{Bou}  Proposition~5  to  the   image   of 
$g(\modx a_{1}\hat{R})$ in $\overline{C}[[T_{n}]]$ we conclude that 
$g(\modx a'_{1}\hat{R})$ is not a  zero  divisor  in  $\hat{R}/a'_{1}\hat{R}$. 
Therefore   $a'_{1},   g$   is   a   regular    $\hat{R}$-sequence.    Since 
$\varphi(a'_{1})=a'_{1}, \varphi(f)=g$, the same holds for $a'_{1},f$ with 
$a'_{1},f\in       \hat{I}$.        We        thus        obtain        that 
$\depth(I,R)=\depth(\hat{I},\hat{R})\geq 2$.
\end{pf}

\begin{eproof}
Recall that we have reduced the proof of the  theorem  to  the  case  of  reduced 
$X=\Spec B$ and we have assumed by way of contradiction  that  
$I=I(Y\setminus 
Y')\subsetneqq R$. We prove that $\depth(I,R)\geq 2$ applying 
Proposition~\ref{s1.8}.  The  condition  $I\not  \subset  \mathfrak{m}R$  is 
fulfilled since  $\mathfrak{m}R=\zeta  \in  Y'$.  We  want  to  verify  that 
Condition~(a) of Proposition~\ref{s1.8} holds. By hypothesis $A$ is reduced, 
so one needs to prove that $A\cap V(I)$ contains none of the minimal prime ideals of 
$A$. Let $X_{i}=V(P_{i}), i=1,\ldots,n$ be  the  irreducible  components  of 
$X$. Assumption~(i) and Assumption~(v) imply, using that $\eta \in  X'$  and 
$\zeta = f(\eta)\in Y'$, that there is a  bijective  correspondence  between 
the irreducible components of $X$ and those of $Y$ given by 
\begin{equation*}
X_{i}\mapsto X_{i}\cap X' \mapsto f(X_{i}\cap X') \mapsto 
\overline{f(X_{i}\cap X')}= f(X_{i}) = Y_{i}.
\end{equation*}
Let us give to $X_{i}\subset  X$  and  $Y_{i}\subset  Y$  the  structure  of 
reduced closed subschemes and let  $f_{i}:X_{i}\to  Y_{i}$  be  the  morphism 
induced by $f$, $i=1,\ldots,n$. For every $i$  the  affine  schemes  $X_{i}, 
Y_{i}$ are integral, the morphism $f_{i}$ is finite and if $x_{i}\in X_{i}$, 
$y_{i}\in Y_{i}$ are the generic points, the homomorphism
\begin{equation*}
(f^{\sharp}_{i})_{y_{i}}: k(y_{i})=\mathcal{O}_{Y_{i},y_{i}} \lto
\mathcal{O}_{X_{i},x_{i}}=k(x_{i})
\end{equation*}
is  an  isomorphism  since  it  coincides   with   $\mathcal{O}_{Y,y_{i}}\to 
\mathcal{O}_{X,x_{i}}$ and $y_{i}\in  Y',  x_{i}=f^{-1}(y_{i})\in  X'$.  Let 
$Y^{\reg}=\{u\in Y|\mathcal{O}_{Y,y}\: \text{is a regular  ring}\}$.  We  claim 
that  $Y^{\reg}\subset  Y'$.  Let  $y\in  Y^{\reg}$.  One   has   that   $y\in 
Y_{i}\setminus  \cup_{j\ne  i}Y_{j}$  for  some  $i$.   The   regular   ring 
$\mathcal{O}_{Y,y}=\mathcal{O}_{Y_{i},y}$ is integrally closed in its  field 
of fractions $k(y_{i})$.  
Let $R_i=\mathcal{O}_{Y_i}(Y_i), B_i=\mathcal{O}_{X_i}(X_i)$, 
$y=\mathfrak{q}\in \Spec R_i$, $S=R_i\setminus \mathfrak{q}$.
The  finite  injective  homomorphism  of  integral 
domains        $f_{i}^{\sharp}(Y_{i}):         \mathcal{O}_{Y_{i}}(Y_{i})\to 
\mathcal{O}_{X_{i}}(X_{i})$  induces  an  isomorphism  of  the   fields   of 
fractions   $k(y_{i})\overset{\sim}{\lto}   k(x_{i})$. Since $S^{-1}R_i$ is integrally 
closed one obtains that $S^{-1}R_i\to S^{-1}B_i$ is an isomorphism,
hence   $f_{i}^{-1}(y)$ consists of a unique point $x\in X_{i}$ and 
$(f^{\sharp}_{i})_{y}: \mathcal{O}_{Y_{i},y}\to \mathcal{O}_{X_{i},x}$ is an 
isomorphism. One has $f^{-1}(y)=f^{-1}_{i}(y)=\{x\}$ since $y\in 
Y_{i}\setminus \cup_{j\ne i}Y_{j}$.  Furthermore  
$\mathcal{O}_{X_{i},x}=\mathcal{O}_{X,x}$ 
since $X$ is reduced and $x\in X_{i}\setminus \cup_{j\ne i}X_{j}$. Therefore 
$(f^{\sharp})_{y}:\mathcal{O}_{Y,y}\to \mathcal{O}_{X,x}$ is an isomorphism, 
so $y\in Y'$. The claim that $Y^{\reg}\subset Y'$ is proved. Suppose now that 
$\mathfrak{q}\in V(I)=Y\setminus Y'$ and $A\cap \mathfrak{q} = \mathfrak{p}$ 
is a minimal prime ideal of $A$. Let $S=A\setminus \mathfrak{p}$. One has by 
Assumption~(b) that 
$R\otimes_{A}k(\mathfrak{p}) = R\otimes_{A}A_{\mathfrak{p}}=S^{-1}R$ is a
regular   ring.   Hence   $R_{\mathfrak{q}}$   is   a   regular   ring, so 
$\mathfrak{q}\in Y^{\reg}$ which contradicts  the  inclusion  $Y^{\reg}\subset 
Y'$ proved above. We thus prove that Condition~(a) of Proposition~\ref{s1.8} 
holds, therefore $\depth(I,R)\geq 2$. Theorem~\ref{s1.1} is proved.
\end{eproof}

\section{A criterion for extending morphisms}\label{s2}
In this section we give the proofs of Theorem~\ref{si2} and Theorem~\ref{a1.1}.  We 
refer to \cite[I. \S 4.4]{D-G} or \cite[II]{SGA1} for the properties of smooth morphisms we  use. 
We use an argument from the course notes of 
M. Mustata of 2009 for the proof of  the following proposition, which we need.
\begin{proposition}\label{v2.1}
Let $g:X\to S$ be a morphism of algebraic varieties  over  an  algebraically 
closed field $k$. Suppose there is an integer $d\geq 1$ such that for  every 
$s\in S$ the fiber $g^{-1}(s)$ is irreducible of dimension $d$.  Then  every 
irreducible component of $S$ is dominated by a unique irreducible  component 
of $X$ and every irreducible component of $X$ is a union of fibers  of  $g$. 
In each of the following three cases there  is  a  bijective  correspondence 
between the irreducible  components  of  $X$  and  those  of  $S$  given  by 
$X_{i}\to \overline{g(X_{i})}$ as well as  by  $S_{i}\to  g^{-1}(S_{i})$  in 
Case~(a).
\renewcommand{\theenumi}{\alph{enumi}}
\begin{enumerate}
\item
$g$ is a closed morphism.
\item
All irreducible components of $X$ have the same dimension.
\item
Every irreducible component of $X$ contains a point in which $g$ is flat.
\end{enumerate}
\end{proposition}
\begin{pf}
\emph{Step  1}.  Let  us  first  suppose  that  $S$  is   irreducible.   Let 
$X=X_{1}\cup \cdots \cup X_{n}$ be an irredundant decomposition of $X$  into 
irreducible  components.  Let  $y\in  S$.  By  hypothesis   $g^{-1}(s)$   is 
irreducible, so  $g^{-1}(s)\subset  X_{j}$  for  some  $j$.  The  hypothesis 
implies that $g$ is  surjective,  consequently  $S=\overline{g(X_{i})}$  for 
some $i$. Renumbering we  may  suppose  that  $i=1$.  Let  $U=X_{1}\setminus 
\cup_{i\geq 2}X_{i}$. If $x\in U$ and $s=g(x)$ then $X_{1}$  is  the  unique 
irreducible component of $X$ which contains $g^{-1}(s)$.  The  image  $g(U)$ 
contains an open subset $V\subset S, V\neq \emptyset$ \cite[Theorem IV.3.7]{Perrin}  and  one  has 
$g^{-1}(s)\subset X_{1}$ for every $s\in V$.  We  claim  that  $X_{i}\subset 
g^{-1}(S\setminus V)$ for $i\geq 2$. Let $x\in X_{i}\setminus  X_{1}$.  Then 
$s=g(x)\notin V$. The inclusion $X_{i}\setminus X_{1}\subset 
g^{-1}(S\setminus V)$ implies $X_{i}\subset g^{-1}(S\setminus V)$. Therefore 
$X_{1}$ is the unique irreducible component of $X$ which dominates $S$.  Let 
$g_{1}=g|_{X_{1}}$. If $s\in V$ then $g_{1}^{-1}(s)=g^{-1}(s)$, so $\dim 
g_{1}^{-1}(s)=d$. This implies that $\dim X_{1}-\dim S=d$ and $\dim 
g_{1}^{-1}(g_{1}(x))\geq d$ for every $x\in X_{1}$ (ibid.). Let 
$x\in X_{1}, s=g(s)=g_{1}(s)$. Then $g_{1}^{-1}(s)$  is  closed  in  $X$  of 
dimension $\geq d$ and moreover it is  contained  in  $g^{-1}(s)$  which  is 
by hypothesis irreducible of dimension $d$. Therefore $g_{1}^{-1}(s)=
g^{-1}(s)$. This shows that $X_{1}$ is a union of fibers of $g:X\to S$. 
\par
\emph{Step 2}. Let now $S$ be arbitrary algebraic  variety.
 Let $X=X_{1}\cup \cdots \cup X_{n}$
$S=S_{1}\cup \cdots \cup  S_{m}$  be  the  irredundant  decompositions  into 
irreducible components. Let $i\in [1,m]$. Since 
$g:X\to S$ is surjective $S_{i}= \overline{g(X_{j})}$ for some $j$.
Applying Step~1 to $S'=S_{i}$ and $X'=g^{-1}(S_{i})$ one concludes that $X_{j}$ 
is the unique irreducible component of $X$ which dominates $S_{i}$. Let 
$X_{i}$ be an arbitrary irreducible component of $X$. Let 
$Z=\overline{g(X_{i})}$. Applying Step~1 to $S'=Z$ and $X'=g^{-1}(Z)$ 
one obtains that $X_{i}$ is a union of fibers of $g:X\to S$.
\par
\emph{Step 3}. \emph{Case~(a)}. If $g:X\to S$ is closed, then 
$S_{i}=g(X_{i})$ 
is closed and irreducible and $X_{i}=g^{-1}(S_{i})$. Hence 
$S=S_{1}\cup \cdots \cup S_{n}$  is  an  irredundant  union  of  irreducible 
closed subsets of $S$. 
\par
\emph{Case (b)}. If $\dim X_{i}=N$ for every $i=1,\ldots,n$ then 
$\dim \overline{g(X_{i})} = N-d$ for every $i=1,\ldots,n$. This implies that 
every $\overline{g(X_{i})}$ is an irreducible component of  $X$.  By  Step~2 
every irreducible component of $S$ is  dominated  by  a  unique  irreducible 
component of $X$, so $S=S_{1}\cup \cdots \cup S_{n}$ with 
$S_{i} = \overline{g(X_{i})}$ is an irredundant decomposition of $S$ into
irreducible components.
\par
\emph{Case (c)}. The set of points $U\subset  X$  in  which  $g$  is  flat  is 
open  \cite[Theorem~53]{Ma1}.  By  hypothesis  $U\cap  X_{i}\neq   \emptyset$   for   every 
$i=1,\ldots,n$. Therefore $U_{i}=U\setminus \cup_{j\neq i}X_{j}$ is an open 
nonempty subset of $X$ contained in $X_{i}$. The  flat  morphisms  are  open 
maps \cite[Theorem~V.5.1]{A-K}, therefore $g(U_{i})$ is an open subset of $S$ contained in the 
irreducible closed subset $\overline{g(X_{i})}$. Therefore 
$S_{i}=\overline{g(X_{i})}$ is an irreducible component of $S$. By Step~2
$S=S_{1}\cup \cdots \cup S_{n}$ with 
$S_{i} = \overline{g(X_{i})}$ is an irredundant decomposition of $S$ into
irreducible components.
\end{pf}
\begin{pot1}
One has to prove that for every $y\in Y$  there  is 
an affine neighborhood of  $y$  such  that  $f$  restricted  to  its  affine 
preimage is an isomorphism. So, it  suffices  to  prove  Theorem~\ref{si2}  for 
affine varieties $X, Y$ and $S$. Let $E = A(X),  D=A(Y),  C=A(S)$.  One  may 
furthermore assume that $\Spec  D\to  \Spec  C$  is  a  smooth  morphism  of 
relative dimension $d$. In order to prove that the finite homomorphism 
$D\to E$ is an isomorphism it suffices to prove that for every maximal ideal 
$\mathfrak{m}_{y}\subset D$, if $T=D\setminus \mathfrak{m}_{y}$, the 
localization $T^{-1}D\to T^{-1}E$ is an isomorphism. Let 
$\mathfrak{m}_{y}\cap C=\mathfrak{m}_{s}, s=h(y)$ and let 
$U=C\setminus \mathfrak{m}_{s}=T\cap C$. Taking into account  that  $E=A(X)$ 
is reduced, we prove that $T^{-1}D\to T^{-1}E$ is an isomorphism applying 
Theorem~\ref{s1.1}, where the local homomorphism $A\to R$ is 
$C_{\mathfrak{m}_{s}}=U^{-1}C\to T^{-1}D=D_{\mathfrak{m}_{y}}$, $B=T^{-1}E$
and the diagram \eqref{es1.1} is 
\begin{equation}
\xymatrix{
\Spec T^{-1}E\ar[rr]^{f}\ar[rd]_{g}&&
\Spec T^{-1}D\ar[dl]^{h}\\
&\Spec U^{-1}C
}
\end{equation}
Let us verify the various assumptions of Theorem~\ref{s1.1}. Assumption~(a)  and 
\linebreak
Assumption~(b) hold since $\Spec D\to \Spec C$ is a smooth morphism and 
$k(\mathfrak{m})=k=k(\mathfrak{n})$. Assumption~(b) of Theorem~\ref{si2}  means 
that $D\to E$ is a  finite  homomorphism  and  $\Spm  E\to  \Spm  D$  is 
surjective. This implies, since $D$ and $E$ are reduced, finitely  generated 
algebras over $k$, that $D\to E$ is injective. By the going-up theorem 
$\Spec E\to \Spec D$ is surjective, hence $\Spec T^{-1}E\to  \Spec  T^{-1}D$ 
is surjective and this shows that Assumption~(i) of Theorem~\ref{s1.1} holds. By 
Assumption~(c) of Theorem~\ref{si2} the fiber $g^{-1}(s)$ is irreducible, hence 
there  is  a  unique  minimal  prime  ideal  $P$  in  $E$   which   contains 
$\mathfrak{m}_{s}E$,  $P=I(g^{-1}(s))$.  The  existence  of  a  point  $x\in 
g^{-1}(s)$ such that $g$ is smooth at $x$ implies that $\Spec E\to \Spec  C$ 
is smooth at $P$, hence the fiber $\Spec E/\mathfrak{m}_{s}E$ is  smooth  at 
 the generic point $P/\mathfrak{m}_{s}E$, i.e. $E_{P}/\mathfrak{m}_{s}E_{P}$ 
is a field. This
implies that  the  ideal  $\mathfrak{m}_{s}E$  has  an  irredundant  primary 
decomposition $\mathfrak{m}_{s}E=Q_{1}\cap Q_{2}\cap  \ldots  \cap  Q_{m}$, 
where $Q_{1}=P$ is a prime ideal and $P\subsetneqq \surd Q_{i}$ for  $i\geq 
2$.  We  claim  that  $P\cap  D\subset  \mathfrak{m}_{y}$.  Indeed,  by  the 
surjectivity of $f:X\to Y$ there exists a point $x\in X$ such that $f(x)=y$, 
so  $\mathfrak{m}_{x}\cap  D\subset   \mathfrak{m}_{y}$.   One   has   $x\in 
g^{-1}(s)$,  so  $\mathfrak{m}_{x}\supset  P$,  therefore  $P\cap   D\subset 
\mathfrak{m}_{y}$. Localizing one obtains a primary decomposition 
\[
\mathfrak{m}_{s}T^{-1}E = T^{-1}P \cap \left(  \cap_{i\geq  2}T^{-1}Q_{i\geq 
2}\right),
\]
where $T^{-1}P\subsetneqq \surd(T^{-1}Q_{i})$ for  $i\geq  2$.  This  proves 
that 
$\Spec T^{-1}E/\mathfrak{m}_{s}T^{-1}E=\Spec B/\mathfrak{m}B$ is irreducible  and 
generically  reduced,  so  Assumption~(ii)  of   \linebreak
Theorem~\ref{s1.1} 
holds. Assumption~(iii) for $\Spec T^{-1}E\to \Spec T^{-1}C$ follows from 
the smoothness of the morphism  $\Spec  E\to  \Spec  C$  at  the  point  $P$ 
mentioned   above.   Assumption~(iv)   follows   from   Assumption~(d)    of 
Theorem~\ref{si2}.  It  remains  to  prove  that  Assumption~(v)   holds.   The 
irreducible components of $\Spec T^{-1}E$ are of the  form  $V(T^{-1}P_{0})$ 
where  $P_{0}$  is  a  minimal  prime  ideal  of  $E$  such  that  $P_{0}\cap 
T=\emptyset$, i.e. $P_{0}\cap D\subset \mathfrak{m}_{y}$.  By  the  going-up 
theorem there is a prime ideal $\mathfrak{p}'$ in E such that 
$P_{0}\subset \mathfrak{p}'$ and 
$\mathfrak{p}'\cap D=\mathfrak{m}_{y}$. Hence $\mathfrak{p}'=\mathfrak{m}_{x}$
for some $x\in X$ with $f(x)=y$. This implies that the irreducible component 
$V(P_{0})$ of $X$ contains a point $x$ of the fiber $g^{-1}(s)$.  The  image 
of the smooth morphism $h:Y\to S$ is open. Let $S_{1}=h(Y)=g(X)$. The fibers 
of  the  morphism  $g:X\to  S_{1}$  are  irreducible  of  dimension  $d$  by 
hypothesis. Using Proposition~\ref{v2.1} we conclude  that  the  irreducible 
component $V(P_{0})$ of $X$, being a union  of  fibers  of  $g:X\to  S_{1}$, 
contains $g^{-1}(s)$, therefore $P_{0}\subset P$. Localizing one obtains 
$T^{-1}P_{0}\subset T^{-1}P$. Therefore $V(T^{-1}P_{0})$ contains 
$V(T^{-1}P)$, which is the preimage of
the  closed  point of the map $\Spec  T^{-1}E\to  \Spec  U^{-1}C$.  All  the 
assumptions of Theorem~\ref{s1.1} were verified, so we conclude that 
$T^{-1}D \to T^{-1}E$ is an isomorphism. Theorem~\ref{si2} is proved.
\end{pot1}
\begin{pot2}
The morphism $f:X\to Y$ is proper and its image contains  $U$,  hence  it  is 
surjective. Furthermore the hypothesis of the theorem implies that  $f$  has 
finite fibers, so $f$  is  a  finite  morphism  by  Zariski's  Main  Theorem 
\cite[Corollary~12.89]{G-W}. One has the commutative diagram \eqref{ei2} with $g=h\circ f : X\to S$. We want to verify that the conditions of Theorem~\ref{si2} hold. Condition~(a) holds by hypothesis. Condition~(b) was verified above. We claim that every fiber $g^{-1}(s), s\in g(X) = h(Y)$ is irreducible of dimension $d$. 
For the proof of this statement we may replace $S$ by its open subset $h(Y)$, so we may assume that $h:Y\to S$ is surjective. The fibers of $h$ are irreducible of dimension $d$, so according to Proposition~\ref{v2.1}(c) the irredundant decompositions of $S$ and $Y$ as finite union of closed irreducible subsets are respectively 
$S=\cup_{i=1}^{m} S_i$ and $Y=\cup_{i=1}^{m} Y_i$, where $S_i = \overline{f(Y_i})$, $\dim Y_i - \dim S_i = d$, and furthermore
every $Y_i$ is a union of fibers of $h:Y\to S$. Let 
$\Gamma_i \subset Y\times_S P$ be the graph of $\varphi|_{Y_i\cap U} : Y_i\cap U \to P$. Then $\Gamma = \cup_{i=1}^m\Gamma_i$ and the irredundant decomposition of $X=\overline{\Gamma}$ is $X = \cup_{i=1}^mX_i$, where $X_i$ is the closure of $\Gamma_i$ in $Y\times_S P$. One has $X\cap (U\times_S P) = \Gamma$, so $f^{-1}(U) = \Gamma$. Let $Z$ be an irreducible component of $g^{-1}(s)$. One has $Z\subset X_i$ for some $i$, hence $s = g(Z) \in S_i$ and $\dim Z \geq \dim X_i - \dim S_i = \dim Y_i - \dim S_i = d$. Since $f:X\to Y$ is finite and $f(Z)\subset h^{-1}(s)$ one has that $f(Z) = h^{-1}(s)$. The map $f|_{\Gamma} : \Gamma \to U$ is an isomorphism, so $Z$ contains the nonempty open subset $f^{-1}(U\cap h^{-1}(s))$, hence $Z$ equals its closure. This shows that $g^{-1}(s)$ is irreducible for every $s\in g(X)=h(Y)$. The remaining conditions of (c) and (d) of Theorem~\ref{si2} hold since $f^{-1}(U) \to U$ is an isomorphism. By  Theorem~\ref{si2} $f: X\to Y$ is an isomorphism. The morphism $\tilde{\varphi} = \pi_{2}\circ f^{-1} : Y\to P$ is the extension of $\varphi : U\to P$. Theorem~\ref{a1.1} is proved.
\end{pot2}

\begin{acknowledgements}
The author is grateful to J. Koll\'{a}r for useful discussions. The author was on leave of absence from the Institute of Mathematics and Informatics of the Bulgarian Academy of Sciences.
\end{acknowledgements}

\bibliography{extmor}

\begin{thebibliography}{10}
\expandafter\ifx\csname url\endcsname\relax
  \def\url#1{\texttt{#1}}\fi
\expandafter\ifx\csname urlprefix\endcsname\relax\def\urlprefix{URL }\fi
\expandafter\ifx\csname href\endcsname\relax
  \def\href#1#2{#2} \def\path#1{#1}\fi

\bibitem{EGAIII.1}
A.~Grothendieck, \'{E}l\'{e}ments de g\'{e}om\'{e}trie alg\'{e}brique. {III}.
  \'{E}tude cohomologique des faisceaux coh\'{e}rents. {I}, Publications
  Math\'ematiques de l'IH\'ES~(11) (1961) 5--167.

\bibitem{Ko}
J.~Koll\'{a}r, Flatness criteria, J. Algebra 175~(2) (1995) 715--727.
\newblock \href {http://dx.doi.org/10.1006/jabr.1995.1209}
  {\path{doi:10.1006/jabr.1995.1209}}.

\bibitem{Ma2}
H.~Matsumura, Commutative ring theory, Vol.~8 of Cambridge Studies in Advanced
  Mathematics, Cambridge University Press, Cambridge, 1986.

\bibitem{Ma1}
H.~Matsumura, Commutative algebra, Vol.~56 of Mathematics Lecture Note Series,
  Benjamin/Cummings, 1980, second edition.

\bibitem{SGA1}
A.~Grothendieck, Rev\^{e}tements \'{e}tales et groupe fondamental (SGA1), Vol.
  224 of Lecture Notes in Mathematics, Springer-Verlag, 1971.

\bibitem{A-K}
A.~Altman, S.~Kleiman, Introduction to {G}rothendieck duality theory, Vol. 146
  of Lecture Notes in Mathematics, Springer-Verlag, 1970.

\bibitem{EGAIV.1}
A.~Grothendieck, \'{E}l\'{e}ments de g\'{e}om\'{e}trie alg\'{e}brique. {IV}.
  \'{E}tude locale des sch\'{e}mas et des morphismes de sch\'{e}mas. {I},
  Publications Math\'ematiques de l'IH\'ES~(20) (1964) 5--259.

\bibitem{Andre}
M.~Andr\'{e}, Localisation de la lissit\'{e} formelle, Manuscripta Math. 13
  (1974) 297--307.
\newblock \href {http://dx.doi.org/10.1007/BF01168230}
  {\path{doi:10.1007/BF01168230}}.

\bibitem{Ei}
D.~Eisenbud, Commutative algebra With a view toward algebraic geometry, Vol.
  150 of Graduate Texts in Mathematics, Springer-Verlag, 1995.
\newblock \href {http://dx.doi.org/10.1007/978-1-4612-5350-1}
  {\path{doi:10.1007/978-1-4612-5350-1}}.

\bibitem{Bou}
N.~Bourbaki, Commutative algebra. {C}hapters 1--7, Elements of Mathematics
  (Berlin), Springer-Verlag, 1998, translated from the French.

\bibitem{A-M}
M.~F. Atiyah, I.~G. Macdonald, Introduction to commutative algebra,
  Addison-Wesley, 1969.

\bibitem{D-G}
M.~Demazure, P.~Gabriel, Groupes alg\'{e}briques. {T}ome {I}:
  {G}\'{e}om\'{e}trie alg\'{e}brique, g\'{e}n\'{e}ralit\'{e}s, groupes
  commutatifs, Masson \& Cie, \'{E}diteur, Paris; North-Holland Publishing Co.,
  Amsterdam, 1970.

\bibitem{Perrin}
D.~Perrin, Algebraic geometry, Universitext, Springer-Verlag London, Ltd.,
  London; EDP Sciences, Les Ulis, 2008, an introduction, Translated from the
  1995 French original by Catriona Maclean.

\bibitem{G-W}
U.~G\"{o}rtz, T.~Wedhorn, Algebraic geometry {I}, Advanced Lectures in
  Mathematics, Vieweg + Teubner, 2010.
\newblock \href {http://dx.doi.org/10.1007/978-3-8348-9722-0}
  {\path{doi:10.1007/978-3-8348-9722-0}}.

\end{thebibliography}

\end{document}